\documentclass[12pt]{article} %% onecolumn (second format) %% use this one for FCAA ! %%

\usepackage{amssymb,amsmath,amsfonts,latexsym} %%
\usepackage{hyperref}

\usepackage{mathtools}
\usepackage{multirow,array}
\usepackage{graphicx}
\usepackage{subfigure}
\usepackage{epstopdf}
\usepackage{float} %%% to use [H], [ht] etc. to fix places of figures
\usepackage{color, xcolor}
\newcommand{\N}{\mathbb{N}}
\newcommand{\R}{\mathbb{R}}
\newcommand{\C}{\mathbb{C}}
\newcommand{\Z}{\mathbb{Z}}

\newcommand{\beq}{\begin{eqnarray}
  }
\newcommand{\eeq}{\end{eqnarray}}
\newcommand{\beqst}{\begin{eqnarray*}
  }
\newcommand{\eeqst}{\end{eqnarray*}}

\newtheorem{theorem}{Theorem}

\newtheorem{lemma}{Lemma}

\newtheorem{remark}{Remark}
\newcommand{\proof}{\noindent{\it Proof}. }
\newcommand{\proofend} {\hfill $\Box$}

\begin{document}

 %%% Pls. specify the kind of the article:
 %%% This is: %% ORIGINAL PAPER %%%% 
 %%% or: REVIEW PAPER / SHORT PAPER, etc. %%%%%%%%

%%%% Title of article for FCAA %%%%%%%%%%%%%%%%%%%%%%%%
\title{Inverse problem to determine simultaneously several scalar parameters and a time-dependent source term in a superdiffusion equation involving a multiterm fractional Laplacian }

\date{}

%%%%% authors:
\author{Hany Gerges, Jaan Janno}

\maketitle

%%%%%%%%%%%%%%%%%%%%%%%%%%%%%%%%%%%%%%%%%%%%%%%%%%%%%%
\begin{abstract}

 Inverse problem to recover simultaneously  a scalar coefficient, order of a time-fractional derivative, parameters of multiterm  fractional Laplacian and a time-dependent source term occurring in a
 superdiffusion equation from  measurements over the time is considered. Uniqueness of a solution is proved. The proof uses asymptotics of poles of the Laplace transform of a measured function.

\end{abstract} %%%%%%%%%%%%%  

%%%%%%%%%% Enter suitable key words and phrases %%% examples:
%%% These are examples only. Pls. use MSC 2020 for suitable topic numbers %%%%%%

%%%%%%%% begin papers' body %%%%%%%%%%%%%%%%%%%%%%%%%%%%%

%%%%%%%%%%%% Section 1 %%%%%%%%%%%%%%%%%%%%%%%%%%
\section{Introduction} \label{sec:1}

%\setcounter{section}{1} \setcounter{equation}{0} %% to have proper 2-digits numbers of eqs
%% Note that this style produces 1-digit numbering of definitons, statements, exmaples, etc.

Differential equations with fractional time derivatives and/or fractional space operators are widely used to model sub- and superdiffusion processes in different media \cite{BKM,Cap,Che,from,spfr1,Sok}. Often
 terms of these equations (parameters, source terms etc.)  are {\it a priori} unknown and cannot be measured directly. However, in many cases they can determined indirectly, making use of measurements of states of processes. 
 This leads to solution of inverse problems. 
 
 Inverse problems to recover parameters of single, multiterm or distributed fractional time derivatives and/or time-dependent source terms by means of time trace measurements 
 are of Volterra type. Existence, uniqueness and stability of solutions can be proved by Titchmarsh convolution theorem, contraction principle or Rothe's method \cite{Ism,J2016,JKas,Lop,Saka,slosis}.
 
 Another class of inverse problems uses space trace measurements (e.g. at a final time). This is convenient in reconstruction of terms depending on space variables. Determination of space-dependent 
 source factors is considered in \cite{Jin,KJ,Kir,orlo,Slo}. The papers \cite{JK,Kian1,Kaif,TU} deal with simultaneous reconstruction of scalar parameters (order of fractional time derivative, unknown final time)  and
 space-dependent terms. Leading asymptotics of Mittag-Leffler functions, Fredholm theory, monotonicity methods etc. are used to prove the uniqueness, stability and in some cases also the existence. 
 
 The paper \cite{JFCAA} is concerned with simultaneous reconstruction of an order of fractional derivative, parameters of distributed fractional Laplacian and a time-dependent source term from final measurements. 
 Full asymptotical expansion of Mittag-Leffler functions is used to prove the uniqueness. 
 
In this paper we will consider a problem to recover simultaneously a scalar coefficient, an order of time-fractional derivative, parameters of multiterm  fractional Laplacian and a time-dependent source term from 
time trace measurements, or more generally, from measurements of a functional of a state over the time. 
We will prove the uniqueness for this problem. The main difficulty is related to the unique recovery of the unknown space operator (i.e. the multiterm fractional Laplacian)
from such measurements. We will use a novel method of asymptotics of poles of the Laplace transform of the measured function. The method is applicable in the case of superdiffusion 
(the order of the time derivative is between $1$ and $2$).  

%%%%%%%%%%%%%%% Section 2 %%%%%%%%%%%%%%%%%%%%%%%
 \section{Problem formulation}\label{sec:form}

%\setcounter{section}{2} \setcounter{equation}{0} %% to have proper 2-digits numbers of eqs
%% Note that this style produces 1-digit numbering of definitons, statements, exmaples, etc.
The Riemann-Liouville fractional integral of the order $\gamma>0$ of a function $w(t)$, $t>0$, is defined by the formula
 $$
 I^\gamma_t w(t)={t^{\gamma-1}\over\Gamma(\gamma)}*w(t),\; t>0,
 $$
 where $*$ denotes the convolution, i.e. $w_1*w_2(t)=\int_0^t w_1(t-\tau)w_2(\tau)d\tau$. This formula admits natural generalization to $\gamma\in\C : {\rm Re}\gamma>0$. Moreover, it is well-known that if $w(t)$ is sufficiently regular 
 then the $\gamma$-dependent function $I_t^\gamma w(t)$ can be analytically continued  to the left complex plane. For $\gamma<0$,
 $I_t^\gamma w(t)$ coincides with the Riemann-Liouville fractional derivative $D_t^\theta w(t)$ of $w$ the order $\theta=-\gamma$. In the particular case $\gamma\in (0,1)$,
  the formula $D_t^\theta w(t)={d\over dt} I_t^{1-\theta}w(t)$ is valid.

Let us consider the time-fractional diffusion equation \cite{BKM,from}
$$
u_t(t,x)=D_t^{1-\alpha}\Delta u(t,x)+f_0(t,x),
$$
where $u$ is a state of a process, $t$ is the time, $x$ is a vector of space variables and $f_0$ is a source function. 
The left-hand side represents the increase rate of a substance at a point $x$ at a time value $t$. The right-hand side represents the
inflow rate to the point $x$ at time $t$, in particular the first addend corresponds to the inflow from neighboring points and time inertia. 
We will consider  the case of superdiffusion (fractional wave process). Then $\alpha\in (1,2)$ and $D_t^{1-\alpha}$ becomes the Riemann-Liouville fractional integral 
 of the positive order $\alpha-1$, i.e. we have
$
u_t(t,x)={t^{\alpha-2}\over\Gamma(\alpha-1)}*\Delta u(t,x)+f_0(t,x).
$
Application of  ${t^{1-\alpha}\over \Gamma(2-\alpha)}*{\partial\over \partial t}$ 
transforms this equation to the equation containing the explicit Laplacian:
\beq\label{Ceq}
\partial_t^{\alpha}u(t,x)=\Delta u(t,x)+{t^{1-\alpha}\over \Gamma(2-\alpha)}*\chi(t,x),\quad \chi=f_{0,t},
\eeq
%%%%parandus%%%
where $\partial_t^{\alpha}u={t^{1-\alpha}\over \Gamma(2-\alpha)}*{\partial^2\over\partial t^2}u$ is the Djrbashian-Caputo derivative of the order $\alpha$ of $u$. 

Let us perform further generalizations of the equation \eqref{Ceq}. To incorporate non-local behavior in space variables, the Laplacian can be replaced by a single, multiterm or distributed fractional elliptic operator
(for such approaches, see \cite{Cap,Che,JFCAA,spfr2,spfr1,Sok}).  We will work with the multiterm spectral fractional Laplacian 
\beq\label{Ldefi}
L=\sum_{j=1}^m b_j (-\Delta)^{\beta_j},
\eeq
where $b_j>0$, $0<\beta_j\le 1$ \cite{JFCAA}. This operator will be defined in the next section. Moreover, note that $\partial_t^{\alpha}u(t,x)=D_t^{\alpha-1}(u_t(t,x)-u_t(0,x))$. Here the right-hand side 
does not involve the second order derivative of $u$, hence applies to more general regular functions $u$. We also introduce a coefficient $a>0$ at the fractional time derivative that may 
occur in non-dimensionless physical models. Summing up, the equation we will consider in this paper is as follows:
%%%%parandus%%%
$
aD_t^{\alpha-1}(u_t-u_t(0,\cdot))+Lu={t^{1-\alpha}\over \Gamma(2-\alpha)}*\chi.
$

Let us proceed to formulation of problems. 
Let $d\in\N$ and $\Omega$ be a bounded domain in $\R^d$ with a sufficiently regular boundary $\partial\Omega$.  Moreover, let $\nu$ denote the outer normal of $\partial\Omega$. 
We introduce the boundary operator ${\mathcal B}$ that acts on functions defined in $\Omega$ and has one of the following two forms:
$$\mbox{${\mathcal B}u(t,x)=u(t,x)$ or ${\mathcal B}u(t,x)={d\over d\nu}u(t,x) +{\rm c} u(t,x)$, \; $t>0$, \, $x\in\partial\Omega$,}
$$
where  ${\rm c}>0$ is a constant. The direct problem consists in determining a function $u$ that satisfies the equation, initial and boundary conditions:
%%%%parandus%%%%%
\beq\label{direc}
\hskip-5truemm \left.
\begin{split}
&aD_t^{\alpha-1}(u_t(t,x)-u_t(0,\cdot))+Lu(t,x)={t^{1-\alpha}\over \Gamma(2-\alpha)}*\chi(t,x),\; x\in\Omega,\, t>0,
\\
&u(0,x)=\varphi(x),\; u_t(0,x)=\psi(x),\; x\in\Omega,
\\
&{\mathcal B}u(t,x)=0,\; x\in\partial\Omega,\, t>0.
\end{split}
\right\}
\eeq

In an inverse problem, we suppose that 
\beq\label{sourcs}
\chi(t,x)=g(t)f(x)+z(t,x),
\eeq
where $g$ is unknown and the scalars $a,\alpha,m,b_j,\beta_j$ are also unknown. We will assume that both $g$ and $z$ have bounded supports in time, i.e. $g(t)=0$, $z(t,\cdot)=0$, $t>T$, where 
$T$ is some positive number. 

To recover the unknowns, we use the observation condition
\beq\label{observs}
\Phi u(t,\cdot )=h(t),\;t\in (0,T+\delta),
\eeq
where $\delta$ is a positive number (it may be arbitrarily small) and $\Phi$ is a functional acting on functions depending on $x$. Examples of $\Phi$ are
\beq\label{Phi1}
&&\Phi w=\int_\Omega \varkappa (x)w(x)dx,\; \varkappa\ge 0,\, \varkappa\not\equiv 0\;\;\mbox{- integral observation}
\\ \label{Phi2}
&&\Phi w=w(x_0),\; x_0\in\Omega\;\;\mbox{- point observation}
\\ \label{Phi3}
&&\Phi w=\int_\Gamma w(x)d\Gamma,\;\Gamma\subseteq\partial\Omega\;\;\mbox{- boundary trace observation}.
\eeq

Summing up, the {\it inverse problem} consists in finding $a,\alpha,m$, $b_1,\ldots,b_m$, $\beta_1,\ldots,\beta_m$, $g$ such that  the solution of the direct problem \eqref{direc} with the source term 
of the form \eqref{sourcs} satisfies the condition \eqref{observs}.

%%%%%%%%  Section 3 %%%%%%%%%%%%%%%%%%%%%%%%%%%%%%
\section{Preliminaries}\label{sec:pre}

%\setcounter{section}{3} \setcounter{equation}{0} %% to have proper 2-digits numbers of eqs
%% Note that this style produces 1-digit numbering of definitons, statements, exmaples, etc.
\subsection{\bf Multiterm fractional powers of minus Laplacian}\label{Lsec}

Let us introduce the eigenvalues and eigenfunctions of $-\Delta$ related to the boundary operator $\mathcal{B}$:
$$
-\Delta v_k(x)=\lambda_k v_k(x),\;\; x\in \Omega,\;\;{\mathcal B}v_k(x)=0,\;\; x\in\partial\Omega.
$$
Let us order $\lambda_k$ in the usual manner: $0<\lambda_1< \lambda_2\le\lambda_3\le \ldots$
and let the system  $(v_k)|_{k\in\N}$ be orthonormed in $L_2(\Omega)$.  It is well-known that $(v_k)|_{k\in\N}$ is complete in $L_2(\Omega)$.

We will use the subscript $k$ to denote Fourier coefficients of functions. This means that
$$
\xi_k=\langle \xi,v_k\rangle_{L_2(\Omega)},\;\; k\in\N,\quad \mbox{for}\quad \xi\in L_2(\Omega).
$$

We define the domain of $-\Delta$ as follows: $\mathcal{D}(-\Delta)=\{\xi\in W_2^2(\Omega)\, :\, {\mathcal B}\xi|_{\partial\Omega}=0\}$. The quantity
$
\|\xi\|_{\mathcal{D}(-\Delta)}=\left[\sum\limits_{k=1}^\infty \lambda_k^2 |\xi_k|^2\right]^{1/ 2}
$
is a norm in $\mathcal{D}(-\Delta)$. 

Let $\rho\ge 0$. The operator $(-\Delta)^\rho$ is defined as follows:
\beq\label{-deltapower}
(-\Delta)^\rho\xi=\sum_{k=1}^\infty \lambda_k^\rho \xi_k v_k.
\eeq
The domain of $(-\Delta)^\rho$ is 
\beq\label{Daste}
\mathcal{D}^\rho:={\mathcal D}((-\Delta)^\rho)=\left\{\xi\in L_2(\Omega)\, :\, \|\xi\|_{{\mathcal D}^\rho}:=\left[\sum_{k=1}^\infty \lambda_k^{2\rho}| \xi_k|^2\right]^{1/2}\!\!\!\!<\infty\right\}\!.
\eeq
This is a Hilbert space with the scalar product
$
\langle \xi,\eta\rangle_{\mathcal{D}^\rho}=\sum\limits_{k=1}^\infty \lambda_k^{2\rho} \xi_k\overline{\eta_k}.
$
It holds 
\beq
\label{emb2}
\mathcal{D}^{\rho_2}\hookrightarrow \mathcal{D}^{\rho_1},\;\;\rho_2>\rho_1\ge 0.
\eeq

Now let us consider the multiterm fractional powers of $-\Delta$. 
Let 
\beq\label{bbeta}
b_j>0,\; j=1,\ldots,m,\quad 1\ge\beta_1>\beta_2>\ldots>\beta_m>0.
\eeq 
Then the operator
$
L=\sum\limits_{j=1}^m b_j (-\Delta)^{\beta_j}
$
is defined in the domain ${\mathcal D}(L)=\bigcap\limits_{j=1}^m \mathcal{D}^{\beta_j}=\mathcal{D}^{\beta_1}$. 
It holds
 \beq\label{Lbounded}
 \mbox{$L$ is a linear bounded operator from $\mathcal{D}^{\rho+\beta_1}$ to $\mathcal{D}^{\rho}$ for $\rho\ge 0$.}
 \eeq

 \subsection{\bf Abstract functional spaces}
Let $X$ be a complex Banach space with the norm $\|\cdot\|$ and $S\subseteq\R$. We define the following spaces of abstract functions 
 mapping $S$ to $X$:
\beqst
&&L_p(S;X)=\left\{w \, :\, \|w\|_{L_p(S;X)}:=\left[\int_S\|w(t)\|^pdt\right]^{1/p}<\infty\right\},\;\; 1\le p<\infty,\\
&&L_\infty(S;X)=\left\{w \, :\, \|w\|_{L_\infty(S;X)}:={\rm ess}\, \sup_{t\in S}\|w(t)\|<\infty\right\},
\\
&&W_p^l(S;X)=\left\{w \, :\, w^{(i)}\in {L_p(S;X)},\, i=0,\ldots,l\right\},\; 1\le p\le \infty,\; l\in\N,
\\
&&C^l(S;X)=\{w\, :\, w^{(i)} - \mbox{continuous in $S$}, i=0,\ldots,l\},\; l\in\N,\\ &&C(S;X)=C^0(S;X).
\eeqst
In case $S$ is compact, $C^l(S;X)$ is a Banach space with the norm $\|w\|_{C^l(S;X)}=\sum\limits_{i=0}^l \max\limits_{t\in S}\|w^{(i)}(t)\|$. 

Next let $T\in (0,\infty)$, $1<p<\infty$, $s>0$ and introduce the following space:
\beqst
&&\hskip-7truemm
_0H_p^s ((0,T);X)=\{w|_{(0,T)}\, :\, w\in H_p^s (\R;X),\, {\rm supp}\,w \subseteq [0,\infty)\},
\eeqst
where $
H_p^s (\R;X)=\{w\in L_p(\R;X)\, :\, {\mathcal F}^{-1}|\xi|^{s}{\mathcal F}w\in L_p(\R;X)\}$
and ${\mathcal F}$ is the Fourier transform with the argument $\xi$.

In case $X=\C$ we drop the value space $\C$ in these notations.

\smallskip
The space $_0H_p^s ((0,T);X)$ is convenient in handling Riemann-Liouville fractional derivatives. 
This can be seen from the following lemma that is a consequence of Corollary 2.8.1 in \cite{zachdiss}.

\begin{lemma}\label{lemma31} 
Let $X$ be a complex Hilbert space.
Let  $\gamma \in (0,1)$, $p\in (1,\infty)$.
The operator  $I_t^{\gamma}$ is a
bijection from $L_p((0,T);X)$ onto $_{0}H_p^\gamma((0,T);X)$, the inverse of $I_t^\gamma$ is the Riemann-Liouville
fractional derivative $D_t^\gamma$
and
$$\|w\|_{_0H_p^\gamma((0,T);X)}=\|{D_t^\gamma } w\|_{L_p((0,T);X)}$$
is a norm in the space $_{0}H_p^\gamma((0,T);X)$.
Moreover,  in case
$p\in ({1\over \gamma},\infty)$ it holds\break\\[-2ex] $_0H_p^\gamma((0,T);X)\hookrightarrow C([0,T];X)$
 and $w(0)=0$ for $w\in {_{0}H_p^\gamma}((0,T);X).$
\end{lemma}

\subsection{\bf Direct problem}\label{sec:dir}

Firstly, we consider the direct problem \eqref{direc} restricted to a bounded time interval $t\in (0,T)$, where  $T\in (0,\infty)$. We introduce the following space for 
solutions of this problem:
\beqst
&&\mathcal{U}_{\alpha,s,\rho,T}=\left\{u\in C^1([0,T];\mathcal{D}^\rho)\cap L_s((0,T);\mathcal{D}^{\rho+\beta_1})\, :\,\right.
\\
&&\qquad \left. u_t-u_t(0,\cdot)\in {_0H_s^{\alpha-1}}((0,T);\mathcal{D}^\rho)\right\},
\quad s>1,\,\rho\ge 0.
\eeqst
In view of \eqref{Lbounded} and Lemma \ref{lemma31}, solutions belonging to $\mathcal{U}_{\alpha,s,\rho,T}$ are strong solutions of \eqref{direc} in $(0,T)\times\Omega$.

\begin{theorem}\label{dpthm}
Let $T\in (0,\infty)$, $a>0$, $\alpha\in (1,2)$, \eqref{bbeta} hold  and consider the direct problem \eqref{direc} restricted to $(t,x)\in (0,T)\times\Omega$.  Then the following assertions are valid.
\begin{enumerate}
\item[\bf (i)]
If the direct problem has a  solution $u\in \mathcal{U}_{\alpha,s,\rho,T}$ for some $s>1$, $\rho\ge 0$, and
$\varphi\equiv 0, \,\psi\equiv 0$, $\chi\equiv 0$ then $u\equiv 0$.
\item[\bf (ii)] If  $\chi|_{(0,T)\times\Omega}$ and $\varphi$, $\psi$ satisfy
\beq\label{chiassum}
&&\hskip-14truemm \chi|_{(0,T)\times\Omega}\in L_p((0,T);\mathcal{D}^{\rho}),\;\sum_{k=1}^\infty \lambda_k^{2\rho}\|\chi_k\|_{L_p(0,T)}^2<\infty,
\\
\label{phipsiassum}
&&\hskip-14truemm\varphi\in \mathcal{D}^{\rho+r_0\beta_1},\, \psi\in \mathcal{D}^{\rho+r_1\beta_1},\, r_j>1\!-\!{j\over \alpha}\!-\!{1\over p\alpha},\, j=0,1,
\eeq 
for some $\rho\ge 0$, $p\in ({1\over\alpha-1},\infty)$,
then the direct problem has a solution $u\in \mathcal{U}_{\alpha,p,\rho,T}$. Fourier coefficients of this solution have the formulas
\beq\label{veel1}
\hskip -5truemm u_k(t)=\varphi_{k}E_{\alpha,1}(-\mu_k t^\alpha)+\psi_{k}tE_{\alpha,2}(-\mu_k t^\alpha)
+{1\over a}[tE_{\alpha,2}(-\mu_k t^\alpha)]\!*\!\chi_k(t), 
\eeq
$k\in\N$, where $E_{\alpha,\vartheta}$ is the two-parametric Mittag-Leffler function and $$\mu_k={1\over a}\sum\limits_{j=1}^m b_l\lambda_k^{\beta_j}.$$
\end{enumerate}
\end{theorem}

\smallskip\noindent Theorem \ref{dpthm} is a slight modification of Theorem 3.1 in \cite{JFCAA}. Its proof  is moved to the Appendix of the paper. 

\smallskip

Let us return to the direct problem \eqref{direc} in the whole time interval $t>0$. 
If $\chi(t,\cdot)=0$, $t>T$, for some $T\in (0,\infty)$  and $\chi|_{(0,T)\times\Omega}$, $\varphi$, $\psi$  satisfy \eqref{chiassum}, \eqref{phipsiassum} for some $\rho\ge 0$, $p\in ({1\over\alpha-1},\infty)$,  then by Theorem  \ref{dpthm}  this problem has a unique solution that satisfies 
$u|_{[0,\tau]\times\Omega}\in \mathcal{U}_{\alpha,p,\rho,\tau}$ for any  $\tau>0$.

\smallskip Let us define the following sectors on the complex plane:
$$
\Sigma(\omega,\theta)=\{z\in\C\, :\, {\rm Re}\, z>\omega,\; |{\rm arg}\, z|< \theta\},\; \omega\in \R,\; \theta\in (0,\pi].
$$

\begin{theorem} \label{dpthm1}
$T\in (0,\infty)$, $a>0$, $\alpha\in (1,2)$, \eqref{bbeta} hold, $\chi(t,\cdot)=0$, $t>T$,  and $\chi|_{(0,T)\times\Omega}$, $\varphi$, $\psi$  satisfy \eqref{chiassum}, \eqref{phipsiassum}  for
some $\rho\ge 0$, $p\in ({1\over\alpha-1},\infty)$. Let $u$ satisfy
$u|_{[0,\tau]\times\Omega}\in \mathcal{U}_{\alpha,p,\rho,\tau}$ for any  $\tau>0$ and solve \eqref{direc}.
Then for any $\sigma>0$ there exists a constant $C_\sigma>0$ such that $e^{-\sigma t}\|u(t,\cdot)\|_{\mathcal{D}^{\rho}}\le C_\sigma$ for $t>0$. Moreover, $u(t,\cdot)$, $t>0$, can be extended to 
a function $u(z,\cdot)$ that is analytic as a function with values in $\mathcal{D}^{\rho}$ in the sector $\Sigma(T,{\pi(2-\alpha)\over 2\alpha}).$
 \end{theorem}

\proof  According to Theorem \ref{dpthm} (ii) and the assumptions of Theorem \ref{dpthm1} the formula \eqref{veel1} is valid for the Fourier coefficients of $u$ for $t>0$.
Since  $E_{\alpha,1}(-t)$ and $E_{\alpha,2}(-t)$ are bounded for $t\ge 0$ and the support of $\chi_k$ is contained in $[0,T]$, for any $\sigma >0$ there exists a constant $C_\sigma^0>0$ such that
\beqst
&&e^{-\sigma t} |\varphi_kE_{\alpha,1}(-\mu_k t^\alpha)|\le C_\sigma^0|\varphi_k|,\; e^{-\sigma t} |\psi_ktE_{\alpha,2}(-\mu_k t^\alpha)|\le C_\sigma^0|\psi_k|,\; 
\\ [1ex]
&&e^{-\sigma t}|a^{-1}[tE_{\alpha,2}(-\mu_k t^\alpha)]*\chi_k(t)|\le C_\sigma^0\|\chi_k\|_{L_p(0,T)},\;\; t>0,\; k\in\N.
\eeqst
Consequently, from \eqref{veel1} we obtain for $\sigma>0$ and $t>0$ that
 \beqst
&&e^{-\sigma t}\|u(t,\cdot)\|_{\mathcal{D}^{\rho}}=e^{-\sigma t}\left[\sum_{k=1}^\infty \lambda_k^{2\rho} |u_{k}(t)|^2\right]^{1\over 2}
\\
&&
\,\le\, C_\sigma^0 \left[\sum_{k=1}^\infty \lambda_k^{2\rho}\left( |\varphi_k|+|\psi_k|+\|\chi_k\|_{L_p(0,T)}\right)^2\right]^{1\over 2}
\\
&&\le C_\sigma^0\left(\|\varphi\|_{\mathcal{D}^{\rho}}+\|\psi\|_{\mathcal{D}^{\rho}}+\left[\sum_{k=1}^\infty \lambda_k^{2\rho}
\|\chi_k\|_{L_p(0,T)}^2\right]^{1\over 2}\right)=:C_\sigma<\infty.
\eeqst

Next let us  consider the functions
\beq\nonumber
&&u_{k}(z)=\varphi_kE_{\alpha,1}(-\mu_k z^\alpha)+\psi_kzE_{\alpha,2}(-\mu_k z^\alpha)
\\ \label{uekz}
&&
+{1\over a}\int_0^{T}(z-\tau)E_{\alpha,2}(-\mu_k (z-\tau)^\alpha)\chi_k(\tau)d\tau, \;\; {\rm Re}\, z>T,\; k\in\N, 
\eeq
where $z^\alpha$ is the principal value of the $\alpha$-th power of $z$. For $z=t$, $t>T$, this function coincides with the 
Fourier coefficient of $u(t,\cdot)$. Since the Mittag-Leffler functions are entire, the functions $E_{\alpha,j}(-\mu_k z^\alpha)$, $j=1,2$,  $k\in\N$, are analytic 
in $\Sigma(0,\pi)$.  Using Lemma 5 of \cite{Janno2024}  we deduce that $$\int_0^{T}(z-\tau)E_{\alpha,2}(-\mu_k (z-\tau)^\alpha)\chi_k(\tau)d\tau$$ is analytic in
$\Sigma(T,{\pi})$. Therefore, $u_{k}(z)$ can be analytically extended to $\Sigma(T,{\pi})$. This implies that 
 the $\mathcal{D}^{\rho}$-valued function
$
\sum\limits_{k=1}^K u_{k}(z)v_k
$
is analytic    in $\Sigma(T,{\pi})$ for any $K\in\N$. 
 
Using the
asymptotic relation for two-parametric Mittag-Leffler functions (6.11) in \cite{Haub}, we deduce that $E_{\alpha,j}(-z)$, $j=1,2$, is bounded in every sector $\Sigma(0,\omega)$, 
$0<\omega<{\pi(2-\alpha)\over 2}$. Therefore, the family of functions $E_{\alpha,j}(-\mu_k z^\alpha)$, $j=1,2$,  $k\in\N$, is uniformly bounded in every sector
$\Sigma(0,\omega)$,  $0<\omega<{\pi(2-\alpha)\over 2\alpha}$. Using this in \eqref{uekz}  we deduce  that for every bounded subset $D$  of $\Sigma(T,\omega)$, where $0<\omega<{\pi(2-\alpha)\over 2\alpha}$,  
 there exists a constant $C_D>0$ such that 
$$|u_{k}(z)|\le C_D\left[|\varphi_k|+|\psi_k|+\|\chi_k\|_{L_p(0,T)}\right],\; z\in D,\; k\in\N.
$$
This yields
\beqst
&&\hskip-5truemm \Big\|\sum_{k=K+1}^{\infty} u_{k}(z)v_k\Big\|_{\mathcal{D}^{\rho}}=\left[\sum_{k=K+1}^{\infty}\lambda_k^{2\rho}  |u_{k}(z)|^2\right]^{1\over 2}
\le C_D\left\{\left[\sum_{k=K+1}^{\infty}\lambda_k^{2\rho}  |\varphi_k|^2\right]^{1\over 2}\right.
\\
&&\hskip-5truemm\left.+
\left[\sum_{k=K+1}^{\infty}\lambda_k^{2\rho}  |\psi_k|^2\right]^{1\over 2}+\left[\sum_{k=K+1}^{\infty}\lambda_k^{2\rho}  \|\chi_k\|^2_{L_p(0,T)}\right]^{1\over 2}\right\}\to 0\;\; \mbox{as}\;\; K\to\infty
\eeqst
for $z\in D$. This shows that $\sum\limits_{k=1}^K u_{k}(z)v_k$ converges uniformly in $D$ to the function $u(z,\cdot):=\sum\limits_{k=1}^\infty u_{k}(z)v_k$. Therefore, $u(z,\cdot)$ is analytic in $D$
as a function with values in $\mathcal{D}^\rho$. Since 
$D$ is  an arbitrary bounded subset of $\Sigma(T,\omega)$ and $\omega\in (0,{\pi(2-\alpha)\over 2\alpha})$ is arbitrary, the function 
$u(z,\cdot)$ is analytic in the sector $\Sigma(T,{\pi(2-\alpha)\over 2\alpha})$. The restriction of $u(z,\cdot)$ to real $z=t$ such that $t>T$ coincides with the solution of 
\eqref{direc}. \proofend

%%%%%%%%%%%%% Section 4 %%%%%%%%%%%%%%%%%%%
\section{Laplace transform of a functional of $u$ and its poles}

In this section we will deduce a formula for the Laplace transform of $\Phi u(t,\cdot)$, where $u$ is the solution of \eqref{direc} and $\Phi\in (\mathcal{D}^{\rho})^*$ for some $\rho\ge 0$, and
describe the set of its poles. We will need this material in the analysis of the inverse problem. 

In the sequel we  will use capitals to denote Laplace transforms of  functions denoted by  lower-case letters. For instance, for $w\in L_{1,loc}((0,\infty);X)$ ($X$ - a complex Banach space)
such that $e^{-\sigma t}w(t)\in  L_1((0,\infty);X)$ for some $\sigma\in\R$, we denote
$$
W(s)=\int_0^\infty e^{-st}w(s)ds,\; {\rm Re}\, s>\sigma.
$$

Moreover, we need to extract a subsequence of non-repeating numbers from $(\lambda_k)|_{k\in\N}$. 
The first eigenvalue $\lambda_1$ is simple, but $\lambda_k$, $k\ge 2$, may in general be multiple. 
Let $k_1=1$ 
and  $k_l$, $l\ge 2$, be the integer such that $\lambda_{k_l+1}>\lambda_{k_l}$ and the vector $\lambda_1,\ldots,\lambda_{k_l}$ contains exactly $l-1$ distinct numbers.
Then  the  sequence $(\lambda_{k_l})|_{l\in\N}$ is strictly increasing and contains all numbers that occur in the set $\{\lambda_s\, :\, s\in\N\}$. In view of $a>0$ and \eqref{bbeta}, 
similar properties are valid for $\mu_{k_l}={1\over a}\sum_{j=1}^m b_j \lambda_{k_l}^{\beta_j}$, too, i.e. $(\mu_{k_l})|_{l\in\N}$ is strictly increasing and contains all numbers that occur in the set $\{\mu_s\, :\, s\in\N\}$.

Additionally, let us denote $\mathcal{K}_l=\{k\in \N\, :\, \lambda_k=\lambda_{k_l}\}$.

\begin{lemma}\label{propo1} 
Let $T\in (0,\infty)$, $a>0$, $\alpha\in (1,2)$, \eqref{bbeta} hold, $\chi(t,\cdot)=0$, $t>T$,  and $\chi|_{(0,T)\times\Omega}$, $\varphi$, $\psi$  satisfy \eqref{chiassum}, \eqref{phipsiassum}  for
some $\rho\ge 0$, $p\in ({1\over\alpha-1},\infty)$. Let $u$ satisfy
$u|_{[0,\tau]\times\Omega}\in \mathcal{U}_{\alpha,p,\rho,\tau}$ for any  $\tau>0$ and solve \eqref{direc}.  Moreover, let 
$\Phi\in (\mathcal{D}^\rho)^*$, $\gamma_k=\Phi v_k$, $k\in \N$, and $h(t)=\Phi u(t,\cdot)$. Assume that 
\beq\label{proa1}
\sum_{k=1}^\infty |\gamma_k| |\varphi_k|<\infty,\; \sum_{k=1}^\infty |\gamma_k| |\psi_k|<\infty,\; \sum_{k=1}^\infty |\gamma_k| \|\chi_k\|_{L_1(0,T)}<\infty.
\eeq
Then $H(s)$ admits a meromorphic extension to the sector $\Sigma(0,\pi)$ and is expressed  by the formula 
\beq\label{proa2}
H(s)={1\over s^{2-\alpha}}\sum_{l=1}^\infty{s\widehat \varphi_{l}+\widehat \psi_{l}+{1\over a}\widehat{\mathcal X}_{l}(s)\over s^\alpha+\mu_{k_l}},
\eeq
where $s^\alpha$ is the main branch of the $\alpha$-the power of $s$, i.e. $s^\alpha=|s|^\alpha e^{i\alpha\, {\rm arg}s}$,  and
\beq\label{proa3}
\widehat \varphi_{l}=\sum_{k\in \mathcal{K}_l}\gamma_k\varphi_k, \;  \widehat \psi_{l}=\sum_{k\in \mathcal{K}_l}\gamma_k\psi_k, \; 
\widehat{\mathcal X}_{l}(s)=\sum_{k\in \mathcal{K}_l}\gamma_k{\mathcal X}_{k}(s).
\eeq
%%%%parandus%%%%%%
Moreover, the extension of $H(s)$ has following set of poles: $\mathcal{P}=\{s_l\in \{\mu_{k_l}^{1\over\alpha}e^{ i{\pi\over\alpha}};\mu_{k_l}^{1\over\alpha}e^{- i{\pi\over\alpha}}\} : s_l\widehat \varphi_{l}+\widehat \psi_{l}+
{1\over a}\widehat{\mathcal X}_{l}(s_l)\ne 0,\, l\in\N\}$. 
\end{lemma}

\proof
Let us consider the expression \eqref{veel1}. Using the formulas of 
Laplace transforms of functions $E_{\alpha,j}(-ct^\alpha)$,  $c>0$, $j=1,2$, (see e.g.  \cite{ML}) we deduce
\beq\label{prop0}
U_k(s)=\int_0^\infty e^{-st}u_{k}(t)dt={1\over s^{2-\alpha}} {s\varphi_{k}+\psi_{k}+{1\over a}{\mathcal X}_{k}(s)\over s^\alpha+\mu_k},\; k\in\N.
\eeq
Moreover, due to Theorem \ref{dpthm1}, $h(t)$ admits a Laplace transform $H(s)$ for ${\rm Re}\, s>0$ and we have
\beq\label{prop1}
H(s)=\int_0^\infty e^{-st}\Phi\sum_{k=1}^\infty u_{k}(t)v_k dt=\int_0^\infty e^{-st}\sum_{k=1}^\infty \gamma_k u_{k}(t)dt,\; {\rm Re}\, s>0.
\eeq
Using the assumptions \eqref{proa1} and the boundedness of $E_{\alpha,j}(-ct^\alpha)$,  $c>0$, $j=1,2$, for $t>0$,  we deduce the estimate
\beqst
&& \int_0^\infty\sum_{k=1}^\infty\left|e^{-st}\gamma_k u_{k}(t)\right|dt
\\
&&\le {\rm const} \int_0^\infty e^{-{\rm Re}s t}\sum_{k=1}^\infty|\gamma_k|\left(|\varphi_k|+t|\psi_k|+\|\chi_k\|_{L_1(0,T)}\right)dt
\\
&&\le {\rm const}  \int_0^\infty e^{-{\rm Re}s t}\max\{1;t\}dt<\infty,\;\; {\rm Re}\, s>0. 
\eeqst
Therefore, due to the Tonelli's theorem we can change the order of  integration and summation in \eqref{prop1}. We obtain 
$H(s)=\sum_{k=1}^\infty \gamma_k \int_0^\infty e^{-st}u_{k}(t)dt$, ${\rm Re}\, s>0$. Using \eqref{prop0} we reach \eqref{proa2} 
for ${\rm Re}\, s>0$. 

Let us choose some $n\in\N$ and $\delta>0$ and denote $$\C_{n,\delta}=\{s\in \Sigma(0,\pi)\, :\, |s^\alpha+\mu_{k_l}|>\delta,\, l\in\N\setminus\{n\}\}.$$ Since $\chi_{k}(t)$, $k\in\N$, 
 have compact supports, ${\mathcal X}_{k}(s)$, $k\in\N$, are entire. Therefore,
for any $l_0\in\N$ the function 
$
\sum_{l=1}^{l_0}{(s\widehat \varphi_{l}+\widehat \psi_{l}+a^{-1}\widehat{\mathcal X}_{l}(s))(s^\alpha+\mu_{k_n})\over s^\alpha+\mu_{k_l}}
$
is analytic in $\C_{n,\delta}$. Moreover, using \eqref{proa1} we deduce that for any $\varrho>0$ and $s\in \C_{n,\delta}$ such that $|s|<\varrho$ the following estimate is valid:
\beqst
&&\left|\sum_{l=l_0+1}^\infty{(s\widehat \varphi_{l}+\widehat \psi_{l}+{1\over a}\widehat{\mathcal X}_{l}(s))(s^\alpha+\mu_{k_n})\over s^\alpha+\mu_{k_l}}\right|
\\
&&\le \max\left\{1;{\varrho^\alpha+\mu_{k_{n}}\over \delta}\right\}\sum_{k=k_{l_0}+1}^\infty |\gamma_k|
\left(\varrho|\varphi_k|+|\psi_k|+{1\over a}\int_0^Te^{{-\rm Re}s t} |\chi_k(t)|dt\right)
\\
&&\le  \max\left\{1;{\varrho^\alpha+\mu_{k_n}\over \delta}\right\}\Bigg(\varrho\sum_{k=k_{l_0}+1}^\infty |\gamma_k||\varphi_k|+\sum_{k=k_{l_0}+1}^\infty |\gamma_k||\psi_k|
\\
&&+
\max\{1;e^{\varrho T}\}{1\over a}\sum_{k=k_{l_0}+1}^\infty |\gamma_k|\|\chi_k\|_{L_1(0,T)}\Bigg)\to 0\;\;\mbox{as}\;\;l_0\to\infty.
\eeqst
This shows that the function $\sum_{l=1}^{l_0}{(s\widehat \varphi_{l}+\widehat \psi_{l}+a^{-1}\widehat{\mathcal X}_{l}(s))(s^\alpha+\mu_{k_n})\over s^\alpha+\mu_{k_l}}$
converges uniformly to the function $$Q_n(s):=\sum_{l=1}^\infty{(s\widehat \varphi_{l}+\widehat \psi_{l}+{1\over a}\widehat{\mathcal X}_{l}(s))(s^\alpha+\mu_{k_n})\over s^\alpha+\mu_{k_l}}$$ in 
every bounded subset of $\C_{n,\delta}$. Therefore, $Q_n(s)$
is analytic in $\C_{n,\delta}$. Since $\delta>0$ is arbitrary, $Q_n(s)$  is analytic 
in the set 
\beqst
&&D_n=\{s\in \Sigma(0,\pi)\, :\, s^\alpha+\mu_{k_l}\ne 0,\, l\in\N\setminus\{n\}\}
\\
&&=\{s\in \Sigma(0,\pi)\, :\, s\ne \mu_{k_l}^{1\over\alpha}e^{\pm i{\pi\over\alpha}},\, l\in\N\setminus\{n\}\}.
\eeqst
Further, let us represent $H(s)$ as
$
H(s)=
{1\over s^\alpha+\mu_{k_n}}Q_n(s)
$
for ${\rm Re}\, s>0$. 
The factor ${1\over s^\alpha+\mu_{k_n}}$ is meromorphic in $\Sigma(0,\pi)$ and has poles at $s_n^\pm= \mu_{k_n}^{1\over\alpha}e^{\pm i{\pi\over\alpha}}$ and $Q_n(s)$ is analytic on $D_n$. 
Consequently, the function $H(s)$ is  meromorphically extendable to 
the set $D_n$ and  $s\in\{s_n^+;s_n^-\}$ is its pole provided  $Q_n(s)\ne 0$. For $s\in\{s_n^+;s_n^-\}$ it holds 
$Q_n(s)\ne 0\Leftrightarrow s\widehat \varphi_{n}+\widehat \psi_{n}+{1\over a}\widehat{\mathcal X}_{n}(s)\ne 0$.  Since $n\in\N$ is arbitrary we reach the assertion of the proposition  about the
meromorphy and poles of $H(s)$. \proofend

\medskip

%%%%%parandus%%%%%
Next we will deduce sufficient conditions that guarantee existence of sufficiently big number of poles of $H(s)$, i.e., the  numbers $s_l\in \{\mu_{k_l}^{1\over\alpha}e^{ i{\pi\over\alpha}}
;\mu_{k_l}^{1\over\alpha}e^{ -i{\pi\over\alpha}}\}$,
$l\in\N$, such that $s_l\widehat \varphi_{l}+\widehat \psi_{l}+
{1\over a}\widehat{\mathcal X}_{l}(s_l)\ne 0$. We will do this assuming $\chi$ has the form \eqref{sourcs} and specify the conditions for  $g,f,z,\varphi,\psi$. Recall that $g$ is unknown in the context of the inverse problem.

In addition to  $\widehat \varphi_l$ and $\widehat \psi_l$ defined in \eqref{proa3}, we introduce the following notation:
\beq\label{proa3a}
\widehat f_l=\sum_{k\in \mathcal{K}_l}\gamma_k f_{k},\;\widehat{z}_{l}=\sum_{k\in \mathcal{K}_l}\gamma_k{ z}_{k},
\eeq
where $\gamma_k=\Phi v_k$, $\Phi\in (\mathcal{D}^\rho)^*$.

We also define the following set of  non-degenerate sequences of real numbers:
$$
\mathcal{ND}=\{(\xi_l)|_{l\in \N}\, :\, \exists A\subseteq \N, \,\mbox{$A$ - infinite,\, $\xi_l\ne 0$, $l\in A$} \}. 
$$

\begin{lemma}\label{invunil1} Let   $T\in (0,\infty)$, $a>0$,  $\alpha \in (1,2)$,   \eqref{bbeta} hold, $\chi$ has the form \eqref{sourcs}, where $g(t)=0$,  $z(t,\cdot)=0$ for $t>T$. Assume that
$g|_{(0,T)}\in L_p(0,T)$, $f\in \mathcal{D}^{\rho}$, 
\beq\label{group3}
z|_{(0,T)\times\Omega}\in L_p((0,T);\mathcal{D}^{\rho}),\;\sum_{k=1}^\infty \lambda_k^{2\rho}\|z_k\|_{L_p(0,T)}^2<\infty
\eeq
 for some $\rho\ge 0$, $p\in ({1\over \alpha-1},\infty)$ and 
 \eqref{phipsiassum} is valid. 
Moreover, let 
$\Phi\in (\mathcal{D}^{\rho})^*$, $\gamma_k=\Phi v_k$ and
\beq\label{uniread}\begin{split}
&\sum_{k=1}^\infty |\gamma_k| |\varphi_k|<\infty,\, \sum_{k=1}^\infty |\gamma_k| |\psi_k|<\infty,\, \sum_{k=1}^\infty |\gamma_k| |f_k|<\infty,\, 
\\
&\sum_{k=1}^\infty |\gamma_k| \|z_k\|_{L_1(0,T)}<\infty.
\end{split}
\eeq
In addition, let  ${\rm Im}\,\varphi={\rm Im}\,\psi=0$, 
%%%parandus%%%%
\beq \label{uniI0}
&&|\widehat \varphi_l|+|\widehat \psi_l|+|\widehat f_l|+\|\widehat z_l\|_{L_1(0,T)}\ne 0,\; l\in\N.
\eeq
Moreover, assume that there exist $l_1\in\N$, $n\in\N$, $C_\dagger>0$  such that for any $l\ge l_1$ the following relations are valid:
\beq
\label{uniz}
&&\widehat z_l|_{[0,T]}\in C^n[0,T],\;\widehat z_l(T)=\ldots=\widehat z_l^{(n-2)}(T)=0 \mbox{ (in case $n\ge 2$)},
\\ \label{uniz0}
&&\|\widehat z_l^{(n)}\|_{C[0,T]}+\sum_{j=0}^{n-1}|\widehat z_l^{(j)}(0)|+|\widehat f_l|\le C_\dagger |\widehat z_l^{(n-1)}(T)|,
\\
\label{eth1}
&&
 \mbox{if\; $\widehat z_l^{(n-1)}(T)\ne 0$\; then\;}\,
|\widehat\varphi_l|+|\widehat \psi_l|\le C_\dagger |\widehat z_l^{(n-1)}(T)|.
 \eeq
In case  $(\widehat f_l)|_{l\in\N}\in\mathcal{ND}$
 we also assume that 
\beq\label{eth2}
g|_{[0,T]}\in C^n[0,T],\; g(T)=\ldots =g^{(n-1)}(T)=0.
\eeq
Then there exists $l_0\in\N$ depending on $a,\alpha,L,g,f,z,\varphi,\psi$ such that 
 for $l\ge l_0$ the following inequality holds:
\beq\label{invunil11} 
s_l\widehat \varphi_{l}+\widehat \psi_{l}+
{1\over a}\left(\widehat f_{l}G(s_l)+\widehat Z_l(s_l)\right)\ne 0,\eeq 
where $s_l=\mu_{k_l}^{1\over\alpha}e^{i{\pi\over\alpha}}$.  Finally, if $l_1=1$, the
conditions \eqref{eth2} hold
and for $l$ such that  $\widehat z_l^{(n-1)}(T)\ne 0$ the following relations are valid: 
%%%%parandus%%%%%
\beq\label{gkits}
\left.\begin{split}
&C_0|\widehat z_l^{(n-1)}(T)|> \|\widehat z_l^{(n)}\|_{C[0,T]}+\sum_{j=0}^{n-1}|\widehat z_l^{(j)}(0)|+|\widehat\varphi_l|+|\widehat \psi_l|,
\\
&|\widehat f_l|\Big(\|g^{(n)}\|_{C[0,T)}+\sum_{j=0}^{n-1}|g^{(j)}(0)|\Big)
\\ 
&<C_0|\widehat z_l^{(n-1)}(T)|- \Big(\|\widehat z_l^{(n)}\|_{C[0,T]}+\sum_{j=0}^{n-1}|\widehat z_l^{(j)}(0)|+|\widehat\varphi_l|+|\widehat \psi_l|\Big),
\end{split}\right\}
\eeq
where
\beq\label{C+form}
C_0= { 
\min\{(\underline b/\overline a)^{1/\underline\alpha};(\underline b/\overline a)^{1/\overline\alpha}\}\min\{\lambda_1^{1/\underline\alpha};
\lambda_1^{\underline\beta/\overline\alpha}\}|\cos(\pi/\overline\alpha)|^{n+2}\over    \max\Big\{1;\left({n+2\over Te}\right)^{n+2}\Big\}\max\{1;\overline a\}},
\eeq 
and $\underline b$, $\overline a$,  $\underline\alpha$, 
$\overline\alpha$, $\underline\beta$, 
 are some numbers satisfying $0<\underline b\le b_j$, $j=1,\ldots,m$, $a\le \overline a$,  $1<\underline\alpha\le\alpha\le \overline\alpha<2$,
$0<\underline\beta\le\beta_j$, $j=1,\ldots,m$, 
then the inequality \eqref{invunil11} holds for any $l\in\N$.

\end{lemma}

\proof Firstly, let us consider the case  $(\widehat f_l)|_{l\in\N}\in\mathcal{ND}$. 
Let $l\ge l_1$. If  $\widehat z_l^{(n-1)}(T)=0$, then by \eqref{uniz0} we have $\hat z_l\equiv 0$ and $\hat f_l=0$, hence
 $s_l\widehat \varphi_{l}+\widehat \psi_{l}+{1\over a}(\widehat f_{l}G(s_l)+\widehat Z_l(s_l))=
s_l\widehat \varphi_{l}+\widehat \psi_{l}$.  Moreover, by \eqref{uniI0}, $|\widehat\varphi_l|+|\widehat\psi_l|\ne 0$. Since $\widehat\varphi_l$ and $\widehat\psi_l$ are real and $s_l$ has nonzero imaginary part, we  
obtain \eqref{invunil11}. 
Next let us suppose that  $\widehat z_l^{(n-1)}(T)\ne 0$.
Due to the property $z(t,\cdot)=0$, $t>T$, we have $\widehat Z_l(s)=\int_0^T e^{-st}\widehat z_l(s)dt$. Integrating by parts $n$ times and taking  \eqref{uniz} into account we deduce
\beqst
&&\widehat Z_l(s)=\sum_{j=0}^{n-1}{1\over s^{j+1}}\widehat z_l^{(j)}(0)-{e^{-sT}\over s^n}\widehat z_l^{(n-1)}(T)+{1\over s^n}\int_0^Te^{-st}\widehat z_l^{(n)}(t)dt.
\eeqst
Let us set $s=s_l$ here. Since ${\rm Re}s_l=\mu_{k_l}^{1\over\alpha}\cos{\pi\over\alpha}<0$  (the factor $\cos{\pi\over\alpha}$ is negative, because $\alpha\in (1,2)$), we deduce
\beq\nonumber
&&\hskip-10truemm  |\widehat Z_l(s_l)|\ge {e^{-{\rm Re}s_lT}\over |s_l|^n}|\widehat z_l^{(n-1)}(T)|- \sum_{j=0}^{n-1}{|\widehat z_l^{(j)}(0)|\over |s_l|^{j+1}}-{e^{-{\rm Re}s_lT}-1\over |s_l|^n|{\rm Re}s_l|}\|\widehat z_l^{(n)}\|_{C[0,T]}
\\ \label{pani1}
&&\hskip-10truemm\ge\! {e^{|{\rm Re}s_l|T}\over |s_l|^n}\!\left(\!|\widehat z_l^{(n-1)}(T)|-{\|\widehat z_l^{(n)}\|_{C[0,T]}\over |{\rm Re}s_l|}-\sum_{j=0}^{n-1}|\widehat z_l^{(j)}(0)| |s_l|^{n-j-1}e^{-|{\rm Re}s_l|T}
\!\right)\!\!.
\eeq
Next by virtue of $g(t)=0$, $t>T$, and \eqref{eth2} we obtain
\beqst
G(s)=\sum_{j=0}^{n-1}{1\over s^{j+1}}g^{(j)}(0)+{1\over s^{n}}\int_0^T e^{-st}g^{(n)}(t)dt.
\eeqst
This implies
\beq\label{pani2}
|G(s_l)|\le {e^{|{\rm Re}s_l|T}\over |s_l|^n} \left({\|g^{(n)}\|_{C[0,T]}\over |{\rm Re}s_l|}+\sum_{j=0}^{n-1}|g^{(j)}(0)| |s_l|^{n-j-1}e^{-|{\rm Re}s_l|T}
\right).
\eeq
By means of \eqref{pani1} and \eqref{pani2} we obtain
\beq\nonumber
&&\left|s_l\widehat \varphi_{l}+
\widehat \psi_{l}+{1\over a}\left(\widehat f_{l}G(s_l)+\widehat Z_l(s_l)\right)\right|\ge {e^{|{\rm Re}s_l|T}\over a|s_l|^n}\Bigg[|\widehat z_l^{(n-1)}(T)|
\\ \nonumber
&&-{\|\widehat z_l^{(n)}\|_{C[0,T]}+
|\widehat f_l|\|g^{(n)}\|_{C[0,T]}\over |{\rm Re}s_l|}
-\sum_{j=0}^{n-1}\left(|\widehat z_l^{(j)}(0)|+|\widehat f_l||g^{(j)}(0)|\right) 
\\\label{tibur1}
&&\times|s_l|^{n-j-1}e^{-|{\rm Re}s_l|T}-a|\widehat\varphi_l| |s_l|^{n+1}e^{-|{\rm Re}s_l|T}-a|\widehat\psi_l| |s_l|^{n}e^{-|{\rm Re}s_l|T}\Bigg].
\eeq
Using  \eqref{uniz0}, \eqref{eth1} we deduce  
\beq\label{tibur1a}
&&\left|s_l\widehat \varphi_{l}+
\widehat \psi_{l}+{1\over a}\left(\widehat f_{l}G(s_l)+\widehat Z_l(s_l)\right)\right|\ge {e^{|{\rm Re}s_l|T}\over a|s_l|^n}|\widehat z_l^{(n-1)}(T)|Q_l,
\\  \nonumber
&&\mbox{where}\;\;Q_l=
1-C_\dagger^1\max_{j=0,\ldots,n+1}\left\{{1\over |{\rm Re}s_l|};|s_l|^j e^{-|{\rm Re}s_l|T}\right\}
\eeq
and $C_\dagger^1$ is a constant independent of $l$.
We note that $|{\rm Re}s_l|=|s_l|\left|\cos{\pi\over\alpha}\right|$ and $|s_l|=\mu_{k_l}^{1/\alpha}\to\infty$ as $l\to\infty$. Therefore, 
 there exists $l_0\ge l_1$
 such that $Q_l>0$ for $l\ge l_0$. Hence
 from \eqref{tibur1a} we conclude that
 \eqref{invunil11} holds for $l\ge l_0$. 

In the case  $(\widehat f_l)|_{l\in\N}\not\in\mathcal{ND}$ there exists $l_2\in\N$ such that $\widehat f_l=0,\, l\ge l_2$. Without restriction of  generality we may assume that $l_1\ge l_2$. Then we repeat the previous computations without the term
$\widehat f_lG$ and obtain  \eqref{invunil11}  for $l\ge l_0$, where $l_0$ is sufficiently big.

It remains to prove the last statement of the lemma. For $l$ such that  $\widehat z_l^{(n-1)}(T)= 0$ we deduce \eqref{invunil11} as %%%parandus%%%% 
in the beginning of the proof. 
So let us consider $l$ such that $\widehat z_l^{(n-1)}(T)\ne 0$ and work with the inequality \eqref{tibur1}. Observing that $|{\rm Re}s_l|=|s_l||\cos{\pi\over\alpha}|$ 
and applying the elementary relations $x^{j}e^{-cx}\le ({j+1\over ce})^{j+1}{1\over x}$, $x>0$, $j\in\N$, $c>0$; $({j+1\over Te})^{j+1}\le \max\{1;({N+1\over Te})^{N+1}\}$, $j=1,\ldots,N$,   
and ${1\over |\cos{\pi\over\alpha}|}\le {1\over |\cos{\pi\over\alpha}|^j}$, $j\in\N$, we obtain
\beq\nonumber
&&\left|s_l\widehat \varphi_{l}+
\widehat \psi_{l}+{1\over a}\left(\widehat f_{l}G(s_l)+\widehat Z_l(s_l)\right)\right|\ge {e^{|{\rm Re}s_l|T}\over a|s_l|^n}\Bigg[|\widehat z_l^{(n-1)}(T)|
\\ \nonumber
&&-{\|\widehat z_l^{(n)}\|_{C[0,T]}+|f_l|\|g^{(n)}\|_{C[0,T]}\over |s_l||\cos{\pi\over\alpha}|^{n+2}}
-{1\over |s_l||\cos{\pi\over\alpha}|^{n+2}}\max\Big\{1;\left({n+2\over Te}\right)^{n+2}\Big\}
\\ \label{tibur2}
&&\times\Bigg(\sum_{j=0}^{n-1}\left(|\widehat z_l^{(j)}(0)|+|f_l||g^{(j)}(0)|\right)+
a|\widehat\varphi_l|+a|\widehat\psi_l| \Bigg)\Bigg].
\eeq
Next we note that $|\cos{\pi\over\alpha}|\ge |\cos(\pi/\overline\alpha)|$ and $|s_l|=\mu_{k_l}^{1\over\alpha}=\left(\sum_{j=1}^m {b_j\over a}\lambda_{k_l}^{\beta_j}\right)^{1\over\alpha}
\ge\break \left(\sum_{j=1}^m {b_j\over a}\lambda_{1}^{\beta_j}\right)^{1\over\alpha}\ge 
\min\{(\underline b/\overline a)^{1/\underline\alpha};(\underline b/\overline a)^{1/\overline\alpha}\}\min\{\lambda_1^{1/\underline\alpha};
\lambda_1^{\underline\beta/\overline\alpha}\}$. Therefore, by \eqref{C+form}
we have
\beqst
&&\hskip-3truemm \left|s_l\widehat \varphi_{l}+
\widehat \psi_{l}+{1\over a}\left(\widehat f_{l}G(s_l)+\widehat Z_l(s_l)\right)\right|\ge {e^{|{\rm Re}s_l|T}\over C_0a|s_l|^n}\Bigg[C_0|\widehat z_l^{(n-1)}(T)|
\\
&&\hskip-3truemm-\|z_l^{(n)}\|_{C[0,T]}-|f_l|\|g^{(n)}\|_{C[0,T]}-
\sum_{j=0}^{n-1}\left(|\widehat z_l^{(j)}(0)|+|f_l||g^{(j)}(0)|\right)-
|\widehat\varphi_l|-|\widehat\psi_l| \Bigg].
\eeqst
Using the assumption \eqref{gkits}  we deduce \eqref{invunil11}.
 \proofend

\section{Uniqueness for  inverse problem}

In this section will formulate and prove a uniqueness theorem for the posed inverse problem. To this end, we need assumptions  that are independent of 
unknown parameters. The condition $p\in ({1\over \alpha-1},\infty)$ is satisfied
provided  $p\in ({1\over \underline\alpha-1},\infty)$ where $\underline\alpha \in (1,2)$ is some number such that $\underline\alpha\le\alpha$. 
 Observing 
\eqref{emb2} we conclude
that sufficient condition that implies   \eqref{phipsiassum} is
\beq\label{group1}
\hskip-1truecm \varphi\in \mathcal{D}^{\rho+r_0},\; \psi\in \mathcal{D}^{\rho+r_1}, \,r_j>1-{j\over \underline\alpha}-{1\over p\underline\alpha},\, j=0,1,
\eeq
for some $\rho\ge 0$, $p\in ({1\over \underline\alpha-1},\infty)$.

\smallskip
We denote the solution of the direct problem \eqref{direc} with $\chi$ of the form \eqref{sourcs} that depends on the quantities $a,\alpha,L,g$ by $u[a,\alpha,L,g](t,x)$.

\begin{theorem}\label{invunithm} 
%%%parandus%%%%%
Let   $T\in (0,\infty)$, $\overline a>0$, $1<\underline\alpha<\overline\alpha<2$, $\underline b>0$, $0<\underline\beta<1$,
\beqst
&&a,\widetilde a\in (0,\overline a],\; \alpha,\widetilde\alpha\in [\underline\alpha,\overline\alpha],\; m,\widetilde m\in \N,
\\
&&b_1,\ldots,b_m\in [\underline b,\infty),\; 1\ge \beta_1>\ldots>\beta_{m}\ge \underline\beta,
\\
&&\widetilde b_1,\ldots,\widetilde b_{\widetilde m}\in [\underline b,\infty),\; 1\ge \widetilde\beta_1>\ldots>\widetilde \beta_{\widetilde m}\ge \underline\beta.
\eeqst
As before, let 
$L=\sum\limits_{j=1}^m b_j (-\Delta)^{\beta_j}$ and define
$\widetilde L=\sum\limits_{j=1}^{\widetilde m} \widetilde b_j (-\Delta)^{\widetilde\beta_j}$. Let 
  $z(t,\cdot)=0$, $t>T$,  $f\in  \mathcal{D}^{\rho}$ and
\eqref{group3}, \eqref{group1} hold
for some $\rho\ge 0$, $p\in ({1\over \underline\alpha-1},\infty)$. Moreover, let 
$\Phi\in (\mathcal{D}^{\rho})^*$, $\gamma_k=\Phi v_k$, 
${\rm Im}\,\varphi={\rm Im}\,\psi=0$,  the conditions  \eqref{uniread}, \eqref{uniI0} hold and 
there  exist $l_1\in\N$, $n\in\N$, $C_\dagger>0$ such that for any $l\ge l_1$ the  relations \eqref{uniz}, \eqref{uniz0}, \eqref{eth1} are valid.
Assume that  $g(t)=0$, $\widetilde g(t)=0$, $t>T$,  and 
\beq\label{eth2thm}
\left.\begin{split}
&g|_{[0,T]}, \widetilde g|_{[0,T]}\in C^n[0,T],\;
\\
&g(T)=\ldots =g^{(n-1)}(T)=\widetilde g(T)=\ldots =\widetilde g^{(n-1)}(T)=0.
\end{split}\right\}
\eeq
In addition, we assume that $(\widehat z_l^{(n-1)}(T))|_{l\in\N}\in \mathcal{ND}$, $\Phi f\ne 0$ and one of the following conditions \eqref{klas}, \eqref{klas1}, \eqref{klas2} is valid: 
 \beq \label{klas} 
&& \left.\begin{split}
&\beta_1-\beta_m<{1\over 2},\;\;\widetilde\beta_1-\widetilde\beta_m<{1\over 2},\;\; 
\\
&\mbox{for any $i\in\N$ there exist $\theta_i>0$, $\zeta_i>0$ and $\bar l_i\in\N$}\;\;
\\
&\mbox{such that}\;\;k_{l+i}^{2\over d}-k_l^{2\over d}\le \theta_i k_l^{1\over d},\; k_{l-i}\ge \zeta_ik_l\;\;\mbox{for}\;\;l\ge \bar l_i;
\end{split}\right\}
\\[2ex]
\label{klas1}
&& 
d=1,\; \Omega=(0,x_1),\; x_1>0,\; \mathcal{B}u=u\;\mbox{on $\partial\Omega$};
\\[2ex]
\label{klas2}
&&  \left.\begin{split}
&l_1=1 \; \mbox{and for $l$ such that  $\widehat z_l^{(n-1)}(T)\ne 0$}
\\
&\mbox{the relations \eqref{gkits} hold for $g$ and $g$ replaced by $\widetilde g$, }
\\
&\mbox{where the constant $C_0$  is given by \eqref{C+form}. }
\end{split}\right\}
 \eeq
Finally, let there exist $\delta>0$ such that 
\beq\label{unieqcond}
\Phi u[a,\alpha,L,g](t)=\Phi u[\widetilde a,\widetilde \alpha,\widetilde L,\widetilde g](t),\; t\in (0,T+\delta).
\eeq
%%%%parandus%%%%%
Then $a=\widetilde a$, $\alpha=\widetilde\alpha$, $m=\widetilde m$, $b_j=\widetilde b_j$, $\beta_j=\widetilde \beta_j$, $j=1,\ldots,m$, and $g=\widetilde g$. 
\end{theorem}

For the proof of this theorem, we need  two additional lemmas. Their proofs can be found in the appendix of the paper. 

\begin{lemma}\label{pohilem0}
There exist constants $c_1>0$ and $c_2\in\R$ (depending on $d$) such that
\beq\label{pohil0}
\lambda_{k_l}=c_1 k_l^{2\over d}\left(1+c_2 k_l^{-{1\over d}}(1+o(1))\right)\;\;\mbox{as $l\to \infty$}.
\eeq
\end{lemma}

\begin{lemma}\label{pohilem1}
Let $m,\widetilde m\in\N$,  ${\rm b}_1,\ldots,{\rm b}_m,{\rm \widetilde b}_1,\ldots, {\rm \widetilde b}_{\widetilde m}>0$, $\beta_1>\ldots>\beta_{m}>0$,
$\widetilde \beta_1>\ldots>\widetilde\beta_{\widetilde m}>0$ and
\beq\label{pohil1}
\sum_{j=1}^{m}  {\rm b}_j\lambda_{k_l}^{\beta_j}+r_{k_l}=\sum_{j=1}^{\widetilde m} {\rm \widetilde b}_j\lambda_{k_l}^{\widetilde\beta_j},\quad l\in\N,
\eeq
where $r_{k_l}=o(\lambda_{k_l}^{\beta_m})$ as $l\to\infty$. Additionally, let 
either $\beta_m\le \widetilde\beta_{\widetilde m}$ or $r_{k_l}=O(1)$ as $l\to\infty$.
Then $m=\widetilde m$, $\beta_j=\widetilde\beta_j$, ${\rm b}_j={\rm \widetilde b}_j$, $j=1,\ldots,m$, and $r_{k_l}=0$, $l\in\N$. 
\end{lemma}

\noindent
{\it Proof of Theorem \ref{invunithm}}.
\,\,  By Theorem \ref{dpthm1}, $\Phi u[a,\alpha,L,g](t)$ and $\Phi u[\widetilde a,\widetilde \alpha,\widetilde L,\widetilde g](t)$ are analytic for $t>T$. Therefore, by analytic continuation, from 
\eqref{unieqcond} we obtain $h(t):=\Phi u[a,\alpha,L,g](t)=\Phi u[\widetilde a,\widetilde \alpha,\widetilde L,\widetilde g](t)$, $t>0$. Using Lemma \ref{propo1}, we deduce
\beq \nonumber
&&H(s)= {1\over s^{2-\alpha}}\sum_{l=1}^\infty{s\widehat \varphi_{l}+\widehat \psi_{l}+{1\over a}\left(\widehat f_{l}G(s)+\widehat Z_l(s)\right)\over s^\alpha+\mu_{k_l}}
\\ \label{unipr1} 
&&= {1\over s^{2-\widetilde\alpha}}\sum_{l=1}^\infty{s\widehat \varphi_{l}+\widehat \psi_{l}+{1\over \widetilde a}\left(\widehat f_{l}\widetilde G(s)+\widehat Z_l(s)\right)\over s^{\widetilde\alpha}+\widetilde \mu_{k_l}},
\eeq
where
$\mu_k=\sum\limits_{j=1}^m {b_j\over a}\lambda_k^{\beta_j}$, as before, and $\widetilde\mu_k=\sum\limits_{j=1}^{\widetilde m} {\widetilde b_j\over \widetilde a}\lambda_k^{\widetilde\beta_j}$. 
In addition, we denote $s_l=\mu_{k_l}^{1\over\alpha}e^{i{\pi\over\alpha}}$, $\widetilde s_l=\widetilde\mu_{k_l}^{1\over\widetilde\alpha}e^{i{\pi\over\widetilde\alpha}}$, $l\in\N$.

Firstly, we prove that $\alpha=\widetilde\alpha$ and 
\beq\label{bbmm}
m=\widetilde m,\; \beta_j=\widetilde\beta_j,\; {b_j\over a}={\widetilde b_j\over \widetilde a},\; j=1,\ldots,m.
\eeq
The simpliest case is \eqref{klas2}. Then, by the last statement of Lemma \ref{invunil1} and Lemma \ref{propo1}, $s_l$, $l\in\N$, are poles of $H(s)$ and $H(s)$ has no more poles in the upper half-plane. 
Similarly, $\widetilde s_l$, $l\in\N$, are poles of $H(s)$ and $H(s)$ has no more poles in the upper half-plane.  Consequently,
$s_l=\widetilde s_l$, $l\in\N$. From the equality of arguments of $s_l$ and $\widetilde s_l$ we get $\alpha=\widetilde \alpha$. Therefore, the equality of modula yields $\mu_{k_l}=\widetilde\mu_{k_l}$, $l\in\N$. Lemma \ref{pohilem1} implies 
\eqref{bbmm}. 

Let us proceed to the cases \eqref{klas} and \eqref{klas1}. Using  Lemmas \ref{invunil1} and \ref{propo1} we deduce that there exists $l_0\in\N$ such that 
 $s_l$, $l\ge l_0$,  are poles of $H(s)$ and $H(s)$ has no more poles 
 in the set  
 %%%%parandus%%%
 $\{s\, :\, |s|\ge |s_{l_0}|,\,{\rm Im}s>0\}$. Let us denote the sequence of all poles of $H(s)$ located in the upper half-plane  by $(p_n)|_{n\in\N}$  and suppose that they ordered so that $|p_{n+1}|>|p_n|$. 
Moreover, let the number of poles of $H(s)$ that are located in the half-disc 
%%%parandus%%%
$\{s\, :\, |s|< |s_{l_0}|,\,{\rm Im}s>0\}$ be $n_0$. Then the subsequences $(p_n)|_{n> n_0}$ and $(s_l)|_{l\ge l_0}$ coincide. This 
implies that there exist $l_*\in\N$, $q_*\in\N$ such that $p_l=s_{l+q_*}$, $l\ge l_*$. 
On the other hand, by Lemmas \ref{invunil1} and  \ref{propo1} there exists $\widetilde l_0\in\N$ such that 
 $\widetilde s_l$, $l\ge \widetilde l_0$,  are poles of $H(s)$ and $H(s)$ has no more poles 
 %%%parandus%%%
 in the set  $\{s\, :\, |s|\ge |\widetilde s_{\widetilde l_0}|,\,{\rm Im}s>0\}$. Arguing as before, we deduce that
 there exist $\widetilde l_*\in\N$, $\widetilde q_*\in\N$ such that $p_l=\widetilde s_{l+\widetilde q_*}$, $l\ge 
\widetilde l_*$. 
Summing up, we have
$
s_{l+q_*}=\widetilde s_{l+\widetilde q_*}$, $ l\ge l_5=\max\{l_*,\widetilde l_*\}$.
From the equality of arguments we obtain $\alpha=\widetilde \alpha$. The equality of modula yields $\mu_{k_{l+q_*}}=\widetilde \mu_{k_{l+\widetilde q_*}}$,  $l\ge l_5$. 
Let us rewrite this in the following form:
\beq\label{poo1}
\sum_{j=1}^m {b_j\over a}\lambda_{k_{l}}^{\beta_j}+r_{k_l}=\sum_{j=1}^{\widetilde m} {\widetilde b_j\over \widetilde a}\lambda_{k_l}^{\widetilde\beta_j},
\;\; l\ge l_6:=l_5+\widetilde q_*,
\eeq
where 
\beq\label{poor}
r_{k_l}=\sum_{j=1}^m {b_j\over a}\lambda_{k_{l+q}}^{\beta_j}-\sum_{j=1}^m {b_j\over a}\lambda_{k_{l}}^{\beta_j}
\eeq
and $q=q_*-\widetilde q_*\in\Z$.

Now we distinguish the cases \eqref{klas} and \eqref{klas1}. Let 
 \eqref{klas} hold. Without loss of generality we assume that 
$\beta_m\le \widetilde \beta_m$ (otherwise we consider the equation $\sum_{j=1}^{\widetilde m} {\widetilde b_j\over \widetilde a}\lambda_{k_l}^{\widetilde\beta_j}+\widetilde r_{k_l}= \sum_{j=1}^m {b_j\over a}\lambda_{k_{l}}^{\beta_j}$ with  $\widetilde r_{k_l}=\sum_{j=1}^{\widetilde m} {\widetilde b_j\over \widetilde a}\lambda_{k_{l-q}}^{\widetilde \beta_j}-\sum_{j=1}^{\widetilde m} {\widetilde b_j\over \widetilde a}\lambda_{k_{l}}^{\widetilde \beta_j}$ instead of \eqref{poo1}).
Using Lemma \ref{pohilem0} we obtain
\beq\label{poo2}
\hskip-5truemm r_{k_l}=\sum_{j=1}^m {b_j\over a}{c_1^{\beta_j}} \left\{[k_{l+q}^{2\over d}+
c_2 k_{l+q}^{{1\over d}}(1+\eta_{l+q})]^{\beta_j}-
[k_l^{2\over d}+c_2 k_{l}^{1\over d}(1+\eta_l)]^{\beta_j}\right\},
\eeq
where $(\eta_l)|_{l\in\N}$ is a sequence of real numbers such that $\eta_l\to 0$
as $l\to\infty$. Using the mean value theorem we have
\beq\label{poo2a}
\hskip-5truemm r_{k_l}=\sum_{j=1}^m {b_j\over a}{c_1^{\beta_j}}\beta_j \xi_{jl}^{\beta_j-1}\left[k_{l+q}^{2\over d}-
k_l^{2\over d}+
c_2 k_{l+q}^{{1\over d}}(1+\eta_{l+q})-c_2 k_{l}^{1\over d}(1+\eta_l)\right],
\eeq
where $\xi_{jl}$ is between $k_l^{2\over d}+c_2 k_{l}^{1\over d}(1+\eta_l)$ and $k_{l+q}^{2\over d}+
c_2 k_{l+q}^{{1\over d}}(1+\eta_{l+q})$. In case $q\ge 0$ we have $k_{l+q}\ge k_l$ and in case $q<0$ by \eqref{klas} it holds $k_{l+q}\ge \zeta_{-q}k_l$, $l\ge \bar l_{-q}$. Hence, there 
exist $C_3>0$ and  $l_7\in \N$ such that 
$\xi_{jl}\ge C_3 k_l^{2\over d}$, $l\ge l_7$, $j=1,\ldots,m$. This implies $ \xi_{jl}^{\beta_j-1}\le C_3^{\beta_j-1}k_l^{{2\over d}(\beta_j-1)}$, $l\ge l_7$, $j=1,\ldots,m$. 
Next let us estimate the terms $k_{l+q}^{2\over d}-k_{l}^{2\over d}$ and ${k^{1\over d}_{l+q}}$ inside the square brackets in \eqref{poo2a}. 
If $q\ge 0$ then by \eqref{klas} we have $0<k_{l+q}^{2\over d}-k_{l}^{2\over d}\le \theta_qk_l^{1\over d}$, $l\ge \bar l_q$, and 
$k_{l+q}^{{1\over d}}=k_l^{1\over d}\sqrt{1+{{k_{l+q}^{{2\over d}}-k_{l}^{{2\over d}}}\over k_l^{2\over d}}}\le 
k_l^{1\over d}\sqrt{1+\theta_qk_{l}^{-{1\over d}}}\le k_l^{1\over d}(1+{1\over 2}\theta_qk_{l}^{-{1\over d}})\le 2k_{l}^{1\over d}$ for sufficiently large $l$. If $q<0$ then have 
$k_{l+q}^{1\over d}\le k_l^{1\over d}$ and
$|k_{l+q}^{2\over d}-k_{l}^{2\over d}|=k_{\hat l-q}^{2\over d}-k_{\hat l}^{2\over d}$, where $\hat l=l+q$. Thus by  \eqref{klas} we get $|k_{l+q}^{2\over d}-k_{l}^{2\over d}|\le \theta_{-q}k_{\hat l}^{1\over d}< \theta_{-q}k_{l}^{1\over d}$, $l\ge \bar l_{-q}-q$. Using the deduced estimates in \eqref{poo2a} we conclude that there exist $C_4>0$ and $l_8\in\N$ such that 
$|r_{k_l}|\le C_4\sum_{j=1}^m k_l^{{1\over d}(2\beta_j-1)}$, $l\ge l_8$. Moreover, due to \eqref{pohil0}, the exists $l_9\in\N$ such that $\lambda_{k_l}\ge {c_1\over 2}k_l^{2\over d}$, $l\ge l_9$. 
Therefore, ${|r_{k_l}|\over \lambda_{k_l}^{\beta_m}}\le C_4\left({2\over c_1}\right)^{\beta_m}\sum_{j=1}^mk_l^{{1\over d}(2\beta_j-2\beta_m-1)}$, $l\ge \max\{l_8;l_9\}$. Due to the assumption 
$\beta_1-\beta_m<{1\over 2}$ in \eqref{klas} and $\beta_j<\beta_1$, $j=2,\ldots,m$, we have $2\beta_j-2\beta_m-1<0$, $j=1,\ldots,m$, and obtain $r_{k_l}=o(\lambda_{k_l}^{\beta_m})$ as $l\to\infty$. 
This means that we can apply Lemma \ref{pohilem1} to \eqref{poo1}. We obtain \eqref{bbmm}.

Next let \eqref{klas1} be valid. Then  $\lambda_k=\left({\pi k\over x_1}\right)^2$, $k\in\N$, and $k_l=l$, $l\in\N$.  
Again, we consider the relation \eqref{poo1} with \eqref{poor}, but this time we assume without restriction of generality that $q\ge 0$ 
(otherwise we consider the equation $\sum_{j=1}^{\widetilde m} {\widetilde b_j\over \widetilde a}\lambda_{k_l}^{\widetilde\beta_j}+\widetilde r_{k_l}= \sum_{j=1}^m {b_j\over a}\lambda_{k_{l}}^{\beta_j}$ with  $\widetilde r_{k_l}=\sum_{j=1}^{\widetilde m} {\widetilde b_j\over \widetilde a}\lambda_{k_{l-q}}^{\widetilde \beta_j}-\sum_{j=1}^{\widetilde m} {\widetilde b_j\over \widetilde a}\lambda_{k_{l}}^{\widetilde \beta_j}$ instead of \eqref{poo1}).
We have $r_{k_l}=\sum_{j=1}^m {b_j\over a} \left({\pi\over x_1}\right)^{2\beta_j}[(l+q)^{2\beta_j}-l^{2\beta_j}]$. By means of the Taylor's formula we obtain
\beq\nonumber
&&r_{k_l}=
\sum_{j=1}^m {b_j\over a} \left({\pi\over x_1}\right)^{2\beta_j}[2\beta_jl^{2\beta_j-1}q+\beta_j(2\beta_j-1)\hat\xi_{jl}^{2\beta_j-2}q^2]
\\ \label{poo3}
&&=\sum_{j=1}^m {b_j\over a} {\pi\over x_1}2\beta_jq\lambda_{k_l}^{\beta_j-{1\over 2}}+
\sum_{j=1}^m {b_j\over a} \left({\pi\over x_1}\right)^{2\beta_j}\beta_j(2\beta_j-1)q^2\hat\xi_{jl}^{2\beta_j-2},
\eeq
where $\hat\xi_{jl}\in (l,l+q)$. The latter relation implies 
 \beq\label{xiasym}
\hskip-5truemm \hat\xi_{jl}^{2\beta_j-2}=l^{2\beta_j-2}(1+o(1))=\left({x_1\over \pi}\right)^{2\beta_j-2} \lambda_{k_{l}}^{\beta_j-1}(1+o(1))\quad \mbox{as\,\, $l\to\infty$}.
\eeq
Let us define $j_0$ so that $\beta_{j_0}>{1\over 2}$ and $\beta_{j_0+1}\le {1\over 2}$. 
(In case $\beta_1\le {1\over 2}$ we take $j_0=0$ and in case $\beta_m>{1\over 2}$ we set $j_0=m$.) Substituting \eqref{poo3} into \eqref{poo1} and using \eqref{xiasym} we obtain
\beq\label{popp1}
\sum_{j=1}^m {b_j\over a}\lambda_{k_{l}}^{\beta_j}+\sum_{j=1}^{j_0}{b_{j-{1\over 2}}\over a}\lambda_{k_{l}}^{\beta_j-{1\over 2}}+
\widehat r_{k_l}=\sum_{j=1}^{\widetilde m} {\widetilde b_j\over \widetilde a}\lambda_{k_l}^{\widetilde\beta_j},\;\;l\ge l_{6},
\eeq
where $b_{j-{1\over 2}}={b_j} {\pi\over x_1}2\beta_jq$ and 
$$
\widehat r_{k_l}=\underbrace{\sum_{j=j_0+1}^{m}{b_j\over a} {\pi\over x_1}2\beta_jq\lambda_{k_{l}}^{\beta_j-{1\over 2}}}_{S_1}+\underbrace{\sum_{j=1}^m {b_j\over a} \left({\pi\over x_1}\right)^{2}\beta_j(2\beta_j-1)q^2 \lambda_{k_{l}}^{\beta_j-1}(1+o(1))}_{S_2}
$$
as $l\to\infty$.\footnote{
Here we assume that any of the sums $\sum_{j=1}^{j_0}$ or $\sum_{j=j_0+1}^{m}$ is absent if its lower limit exceeds the upper one.} In 
\eqref{popp1} the  coefficients ${b_{j-{1\over 2}}\over a}$, $j=1,\ldots,j_0$, are nonnegative and the exponents 
$\beta_{j}-{1\over 2}$, $j=1,\ldots,j_0,$ belong to $ (0,1]$. 
Moreover, $\widehat r_{k_l}=O(1)$ as $l\to\infty$, because $\beta_{j}-{1\over 2}\le 0$, $j=j_0+1,\ldots,m$, and $\beta_{j}-1\le 0$, $j=1,\ldots,m$.  Therefore, 
after proper rearrangement of addends under the sums $\sum_{j=1}^m$ and $\sum_{j=1}^{j_0}$ in  \eqref{popp1} it satisfies the assumptions of Lemma 
\ref{pohilem1}. We obtain  $\widehat r_{kl}=0$, $l\ge l_{6}$. Let us analyze this relation and show that it implies $q=0$. 
If $\beta_{1}\le {1\over 2}$ then  $j_0=0$ and the first addend of $S_1$ dominates in $\widehat r_{kl}$ in the process $l\to\infty$. Therefore, due to $\widehat r_{kl}=0$ we have
${b_1\over a} {\pi\over x_1}2\beta_1q=0$, thus $q=0$. If $\beta_1>{1\over 2}$ then one of the following 4 cases occurs: (a): $\beta_m>{1\over 2}$; (b): $\beta_1>\beta_{j_0+1}+{1\over 2}$, 
(c): $\beta_1<\beta_{j_0+1}+{1\over 2}$ and (d): $\beta_1=\beta_{j_0+1}+{1\over 2}$. In case (a) the term $S_1$ is absent and the first addend of $S_2$ dominates. We get 
${b_1\over a} \left({\pi\over x_1}\right)^{2}\beta_1(2\beta_1-1)q^2=0$, hence $q=0$. In case (b) also the first addend of $S_2$ dominates and similarly we deduce $q=0$. 
In case (c) the first addend of $S_1$ dominates and we obtain ${b_{j_0+1}\over a} {\pi\over x_1}2\beta_{j_0+1}q=0\Rightarrow q=0$. Finally in  case (d) the first addends of $S_1$ and $S_2$ are of the same order
and we deduce ${b_{j_0+1}\over a} {\pi\over x_1}2\beta_{j_0+1}q+{b_1\over a} \left({\pi\over x_1}\right)^{2}\beta_1(2\beta_1-1)q^2=0$. Left hand side of this equality contains  positive addends 
if $q>0$. Therefore, $q=0$. Summing up, we have shown that $q=0$. Let us return to \eqref{poo1}. Since $q=0$, from \eqref{poor} we have $r_{kl}=0$. Applying Lemma \ref{pohilem1} we obtain
\eqref{bbmm}. 

%%%parandus%%%
The relations $\alpha=\widetilde\alpha$ and \eqref{bbmm} have been proved in all cases. 

Secondly, let us prove that $a=\widetilde a$. The proved relations $\alpha=\widetilde\alpha$ and \eqref{bbmm} imply that $\mu_{k_l}=\widetilde \mu_{k_l}$ and $s_l=\widetilde s_l$ for $l\in\N$. 
Therefore, from \eqref{unipr1} we have 
\beq\label{unipr1*}
\sum_{l=1}^\infty{\left({1\over a}-{1\over \widetilde a}\right)\widehat Z_l(s)+  \widehat f_{l}\left({1\over  a}G(s)-{1\over \widetilde a}\widetilde G(s)\right)\over s^\alpha+\mu_{k_l}}=0.
\eeq
Choosing an arbitrary $l_*\in\N$ and sending $s$ to $s_{l_*}$ we obtain 
the relation 
\beq\label{unipr1*1}
\left({1\over a}-{1\over \widetilde a}\right)\widehat Z_{l}(s_{l})+ 
 \widehat f_{l}\left({1\over  a}G(s_{l})-{1\over \widetilde a}\widetilde G(s_{l})\right)=0
\eeq
for $l=l_*$. Since $l_*$ is arbitrary, \eqref{unipr1*1} holds for any $l\in\N$. Since $(\widehat z_l^{(n-1)}(T))|_{l\in\N}\break \in \mathcal{ND}$, there exists an infinite subset $A\subseteq\N$ such that $\widehat z_{l}^{(n-1)}(T)\ne 0$, $l\in A$.  We use the estimates \eqref{pani1} for $\widehat Z_l(s_l)$ and \eqref{pani2} for $G(s_l)$ and  $\widetilde G(s_l)$ in case 
$l\in A$ as well as the relation $ |\widehat f_l|\le C_\dagger |\widehat z_l^{(n-1)}(T)|$ following from \eqref{uniz0}. These estimates show that there exist $l_{10}\in\N$ and  $C>0$ such that
$|\widehat Z_{l}(s_{l})|\ge {e^{-{\rm Re}s_{l}T}\over 2|s_{l}|^n}|\widehat z_{l}^{(n-1)}(T)|>0$ and $|\widehat f_{l}| |G(s_{l})|
\le {Ce^{-{\rm Re}s_{l}T}\over |s_{l}|^n|{\rm Re} s_{l}|}|\widehat z_{l}^{(n-1)}(T)| $, 
$|\widehat f_{l}| |\widetilde G(s_{l})|\le {Ce^{-{\rm Re}s_{l}T}\over |s_{l}|^n|{\rm Re} s_l|}|\widehat z_{l}^{(n-1)}(T)| $
for $l\ge l_{10}$, $l\in A$. Suppose that $a\ne \widetilde a$. Then, setting $l\in A$ in \eqref{unipr1*1} and letting $l\to \infty$, the first addend on the left hand side differs from zero and dominates over the second one.
This leads to the contradiction. Consequently, $a=\widetilde a$. 

It remains to prove that $g=\widetilde g$. Using $a=\widetilde a$ in 
\eqref{unipr1*} we deduce $(G(s)-\widetilde G(s))\sum\limits_{l=1}^\infty {\widehat  f_{l}s^\alpha\over s^\alpha+\mu_{k_l}}=0$. We have
 $\sum\limits_{l=1}^\infty {\widehat  f_{l}s^\alpha\over s^\alpha+\mu_{k_l}}=\sum\limits_{k=1}^\infty \gamma_k{ f_{k}s^\alpha\over s^\alpha+\mu_{k}}\to \sum\limits_{k=1}^\infty \gamma_k f_{k}=\Phi f$ as $s\in\R$, $s\to\infty$. Therefore, since
 $\Phi f\ne 0$, there exists $s_+\in \R$ such that 
$G(s)-\widetilde G(s)=0$ for $s>s_+$. By analytic continuation we have $G(s)=\widetilde G(s)$, $s\in\C$. This implies $g=\widetilde g$. 
Proof is complete. \proofend

\begin{remark}
%%%%parandus%%%%%
For the observation functionals  defined in Section \ref{sec:form} we have the following relations:
\begin{itemize}
\item functional \eqref{Phi1}: 
$\Phi\in (\mathcal{D}^0)^*=(L_2(\Omega))^*$ if $\varkappa\in L_2(\Omega)$;
\item functional \eqref{Phi2}:
$\Phi\in (\mathcal{D}^\rho)^*$, $\rho={1\over 2}\lfloor{{d\over 2}+1}\rfloor$;
\item functional \eqref{Phi3}: $\Phi\in (\mathcal{D}^{1/4})^*$. 
\end{itemize}
For \eqref{Phi2}, \eqref{Phi3} this follows from the facts that $\mathcal{D}^{1/2}$ and $\mathcal{D}^{1/4}$ are subspaces of $H^1(\Omega)$ and $H^{1/2}(\Omega)$, respectively \cite{Saka,Lischke}, 
%%%parandus%%%
the solution of the problem $\left\{\!\!\begin{array}{ll}\Delta w=\omega\; &\mbox{in $\Omega$},
\\ \mathcal{B}w=0\; &\mbox{on $\partial\Omega$,}\end{array}\right.$ belongs to $H^{k+2}(\Omega)$ if $\omega\in H^k(\Omega)$ \cite{Grisvard}  and embedding theorems.
\end{remark}

\begin{remark}
Although there are some well-known  domains $\Omega$, where the  Laplacian has multiple eigenvalues (e.g. disks and balls), for "most" domains the eigenvalues $\lambda_k$ are simple \cite{uhl}. In the
case of simple eigenvalues it holds
$k_l=l$  and the conditions for $k_l$ in \eqref{klas}  are satisfied with the parameters $\bar l_i=2i$, 
$\theta_i={2i\over d}({3\over 2})^{{2\over d}-1}$, $\zeta_i={1\over 2}$. 
\end{remark}

\begin{remark} The proof  of Theorem \ref{invunithm}  in the case \eqref{klas} uses the relation \eqref{pohil0} that follows from the second order asymptotics of the eigenvalues 
deduced in \cite{ivrii}. An improvement of such an asymptotics, namely a specification of a higher order terms may enable to relax the restrictions $\beta_1-\beta_m<{1\over 2}$, $\widetilde\beta_1-\widetilde\beta_m<{1\over 2}$ in \eqref{klas}. However, these restrictions may be relevant in some physical situations. For instance, when a model with a single fractional Laplacian 
$L=(-\Delta)^\beta$ needs a slight improvement, then one can replace it by a model that involves the operator $L=\sum\limits_{j=1}^m b_j (-\Delta)^{\beta_j}$, where 
$\beta_j\in (\beta-\varepsilon,\beta+\varepsilon)$ and $\varepsilon$ is a small positive number. 
\end{remark}

\begin{remark} Results of the paper can be generalized to the case when the assumption \eqref{uniI0} does not hold, i.e. when $I_0:=\{l\,:\, |\widehat \varphi_l|+|\widehat \psi_l|+|\widehat f_l|+\|\widehat z_l\|_{L_1(0,T)}=0\}\ne \emptyset$. Then one can work in the basic space\break  $X=\overline{{\rm span}\{v_k\, :\, k\in\mathcal{K}_l,\, l\in \N\setminus I_0\}}$ instead of $L_2(\Omega)$ and define 
proper subspaces of $\mathcal{D}^\rho$ there. In \eqref{klas} the whole sequence of $(k_l)_{l\in\N}$ must be replaced by its subsequence, where the elements $k_l$, $l\in I_0$, are removed and in 
\eqref{klas1} the restriction that $I_0=\{n_0-1+n_1i\, :\, i\in\N\}$ with some $n_0,n_1\in\N$ must be added. 
\end{remark}

\begin{remark}
The proof of Theorem \ref{invunithm} cannot be directly generalized to the case  $\alpha,\widetilde\alpha\in (0,1)$, because then the extension of $H(s)$ is analytic in the 
sector $\Sigma(0,\pi)$, i.e. has no poles there. 
\end{remark}

%%%%%%%%%%   Example for Appendices %%%%%%%%%%%%%%%%%%%%%%%%%%
%%%%%%%%%%   but please try to avoid, unless the article's structure needs so %%%%%%%%%%%

\section{Appendix}

\noindent {\it Proof of Theorem \ref{dpthm}}\,\,  The assertion (i) and the assertion (ii) in the case $\chi\equiv 0$ follow from Theorem 3.1 in J. Janno, FCAA, 2020. So, it remains to prove (ii) 
in the case $\varphi\equiv 0$, $\psi\equiv 0$. Let $k\in\N$ and consider the problem
\beq
\label{dpf1}
&&D_t^{\alpha-1}w'_k(t)+\mu_k w_k(t)={1\over a}{t^{1-\alpha}\over \Gamma(2-\alpha)}*\chi_k(t),\quad t>0,
\\
\label{dpf2}
&&w_k(0)=0,\quad w_k'(0)=0.
\eeq
It is well-known that $\omega_k(t)=[t^{\alpha-1}E_{\alpha,\alpha}(-\mu_k t^\alpha)]*\eta(t)$ solves the equation $D_t^{\alpha-1}\omega'_k(t)+\mu_k \omega_k(t)=\eta(t)$, $t>0$. 
Since $[t^{\alpha-1}E_{\alpha,\alpha}(-\mu_k t^\alpha)]*{t^{1-\alpha}\over \Gamma(2-\alpha)}=tE_{\alpha,2}(-\mu_k t^\alpha)$ \cite{ML}, we see that the function 
\beq\label{solutt}
w_k(t)={1\over a}[tE_{\alpha,2}(-\mu_k t^\alpha)]*\chi_k(t)
\eeq
is a solution of \eqref{dpf1}. One can immediately check that \eqref{solutt} satisfies also the initial conditions \eqref{dpf2}. 

Let us estimate \eqref{solutt}. To this end we use the inequality
$|E_{\alpha,2}(-t)|\le {c\over 1+t}$, $t>0$,
 that holds with some constant $c$ (see \cite{JFCAA}), the relation $\lambda_k^{\beta_1}\le {a\over b_1}\mu_k$, $k\in\N$,  and Young's inequality for convolutions. We obtain
\beq \nonumber
&&\hskip-1truecm\lambda_k^{\rho+\beta_1} |w_k(t)|\le{ca\over b_1}\lambda_k^{\rho} \int_0^t {\mu_k(t-\tau)^\alpha
\over 1+\mu_k (t-\tau)^\alpha}(t-\tau)^{1-\alpha}
|\chi_k(\tau)|d\tau
\\ \label{veel3}
&&\hskip-1truecm\quad\le
{ca\over b_1}\lambda_k^{\rho}\int_0^t
(t-\tau)^{1-\alpha}
|\chi_k(\tau)|d\tau\le {ca\over b_1}\lambda_k^{\rho} \|t^{1-\alpha}\|_{L_1(0,T)}\|\chi_k\|_{L_p(0,T)}.
\eeq

Now let us consider the function $w=\sum_{k=1}^\infty w_kv_k$. Using   \eqref{chiassum} and \eqref{veel3} we obtain
\beqst
&&\|w\|_{L_p((0,T);\mathcal{D}^{\rho+\beta_1})}\le 
{ca\over b_1} \|t^{1-\alpha}\|_{L_1(0,T)} \left[\sum_{k=1}^\infty  \lambda_k^{2\rho}\|\chi_k\|_{L_p(0,T)}^2\right]^{1\over 2}<\infty.
\eeqst
Therefore, $w\in L_p((0,T);\mathcal{D}^{\rho+\beta_1})$. The latter relation together with \eqref{dpf1} and \eqref{dpf2} implies that $w$ belongs to $ \mathcal{U}_{\alpha,p,\rho,T}$ and  solves the direct 
problem in case $\varphi\equiv\psi\equiv0$. The proof of such an implication is almost identical with the one included Theorem 3.1 in \cite{JFCAA}. Therefore, we will not repeat it here. 
Setting $u=w$, the existence assertion is proved. We also have $u_k=w_k$. 
Complementing \eqref{solutt} by the term $\varphi_{k}E_{\alpha,1}(-\mu_k t^\alpha)+\psi_{k}tE_{\alpha,2}(-\mu_k t^\alpha)$ that occurs in the case of 
non-vanishing $\varphi$ and $\psi$ (see again Theorem 3.1 in \cite{JFCAA}) we deduce \eqref{veel1}. \proofend

\medskip\noindent
{\it Proof of Lemma \ref{pohilem0}}\,\,
Let $N(\lambda)$ denote the number of eigenvalues $\lambda_k$  (counting multiplicities) that are less than or equal to $\lambda$. It holds \cite{ivrii}
\beq\label{pohil0*}
N(\lambda)=c_1^*\lambda^{d\over 2}\left(1+c_2^*\lambda^{-{1\over 2}}(1+o(1))\right)\;\;\mbox{as $\lambda\to \infty$},
\eeq
where $c_1^*>0$ and $c_2^*\in\R$ are some constants depending on $d$. By the definition of $k_l$, this immediately implies
\beqst
k_l=c_1^*\lambda_{k_l}^{d\over 2}\left(1+c_2^*\lambda_{k_l}^{-{1\over 2}}(1+o(1))\right)\;\;\mbox{as $l\to \infty$}.
\eeqst
Therefore, $\lambda_{k_l}^{d\over 2}=(c_1^*)^{-1}k_l\left(1+c_2^*\lambda_{k_l}^{-{1\over 2}}(1+o(1))\right)^{-1}$ as $l\to\infty$. 
We obtain
\beqst
&&\lambda_{k_l}=(c_1^*)^{-{2\over d}}k_l^{2\over d}\left(1+c_2^*\lambda_{k_l}^{-{1\over 2}}(1+o(1))\right)^{-{2\over d}}
\\
&&=(c_1^*)^{-{2\over d}}k_l^{2\over d}\left(1-{2\over d}c_2^*\lambda_{k_l}^{-{1\over 2}}(1+o(1))\right)
\\
&&=(c_1^*)^{-{2\over d}}k_l^{2\over d}\left(1-{2\over d}c_2^*\left[(c_1^*)^{-{2\over d}}k_l^{2\over d}\left(1-{2\over d}c_2^*\lambda_{k_l}^{-{1\over 2}}(1+o(1))\right)\right]^{-{1\over 2}}(1+o(1))\right)
\\
&&=(c_1^*)^{-{2\over d}}k_l^{2\over d}\left(1-{2\over d}c_2^*(c_1^*)^{{1\over d}}k_l^{-{1\over d}}\left[1-{2\over d}c_2^*\lambda_{k_l}^{-{1\over 2}}(1+o(1))\right]^{-{1\over 2}}(1+o(1))\right)
\\
&&=(c_1^*)^{-{2\over d}}k_l^{2\over d}\left(1-{2\over d}c_2^*(c_1^*)^{{1\over d}}k_l^{-{1\over d}}\left[1+{1\over d}c_2^*\lambda_{k_l}^{-{1\over 2}}(1+o(1))\right](1+o(1))\right)
\eeqst
as $l\to\infty$. This implies \eqref{pohil0}.  \proofend

\medskip\noindent
{\it Proof of Lemma \ref{pohilem1}}\,\,
The terms  ${\rm b}_1\lambda_{k_l}^{\beta_1}$ and ${\rm \widetilde b}_1\lambda_{k_l}^{\widetilde\beta_1}$   dominate on the left- and right-hand sides of  \eqref{pohil1}, respectively,  
 in the process $l\to\infty$. This implies $\beta_1=\widetilde\beta_1$ and ${\rm b}_1={\rm \widetilde b}_1$. The relation  \eqref{pohil1} reduces to
 \beqst
\sum_{j=2}^{m} {\rm b}_j\lambda_{k_l}^{\beta_j}+r_{k_l}=\sum_{j=2}^{\widetilde m} {\rm\widetilde b}_j\lambda_{k_l}^{\widetilde\beta_j},\quad l\in\N.
\eeqst
Now the dominating terms  on the left- and right-hand sides are  ${\rm b}_2\lambda_{k_l}^{\beta_2}$   and ${\rm\widetilde b}_2\lambda_{k_l}^{\widetilde\beta_2}$, respectively. 
Therefore, $\beta_2=\widetilde\beta_2$ and ${\rm b}_2={\rm \widetilde b}_2$.
Repeating this procedure we obtain
$\beta_j=\widetilde\beta_j$, ${\rm b}_j={\rm \widetilde b}_j$, $j=1,\ldots,\min\{m;\widetilde m\}$ and
\beq \label{pohil3}
&&\sum_{j=\widetilde m+1}^{m}  {\rm b}_j\lambda_{k_l}^{\beta_j}+r_{k_l}=0,\quad l\in\N,\quad\mbox{
in case $m>\widetilde m$}, 
\\  \label{pohil4}
&&
r_{k_l}=\sum_{j= m+1}^{\widetilde m} {\rm \widetilde b}_j\lambda_{k_l}^{\widetilde\beta_j},\quad l\in\N,\quad\mbox{
in case $m<\widetilde m$}, 
\\ \label{pohil5}
&& r_{k_l}=0,
\quad l\in\N,\quad\mbox{
in case $m=\widetilde m$}.
\eeq
The cases \eqref{pohil3} and \eqref{pohil4} are contradictive, because the left- and right-hand sides are of different order as $l\to\infty.$ Therefore, only \eqref{pohil5} remains. This 
completes the proof. \proofend

\vskip 5truemm

\noindent{\bf Acknowledgement}.\\[1ex]
  The study was supported by Estonian Research Council grant  PRG832.

%%%%%%%%%%%%%%%%%% REFERENCES: %%%%%%%%%%%%%%%%%%%%%%%%%%%%%%%%%%%%%%%%%%%%
%% BibTeX users: please use \bibliographystyle{spmpsci} %% for math. and phys. sci.
%% Non-BibTeX users: please use the model as below !!! %%

%%%% for FCAA - pls. include directly the Refs items here ! %%%
%%%% following STRICTLY the models below %%%%%
%%%% and ARRANGE the items in ALPHABETIC ORDER for authors' family names !!!

 %%%%%%%%%%%%%%%%%%%%%%%%%%%%%%

\end{document}